\input amstex
\documentstyle{amsppt}
\NoBlackBoxes
\magnification 1200
\hsize 6.5truein
\vsize 8truein
\topmatter
\title On Sums, Products , and the multidimensional Falconer Problem \endtitle
\author  Nets Hawk Katz, \endauthor
\affil Indiana University \endaffil
 \thanks  The  author was supported by NSF grant DMS 0432237
\endthanks 
\endtopmatter

\head \S 0 Introduction \endhead

The Falconer distance conjecture says that if a set $E \subset \Bbb R^d$ has Hausdorff dimension
at least ${d \over 2}$ then its distance set $D=\{|x-y|, x,y \in E\}$ has Hausdorff dimension at least 1.
This problem has received a great deal of attention in recent years. Wolff [W] got a partial result for
$d=2$ showing the conjecture holds if the dimension of the set is at least ${4 \over 3}$. Erdogan
[E] extended Wolff's method to obtain partial results in higher dimensions. 

In this paper, we prove essentially

\bigskip

{\bf Theorem 0.1} There exists $\epsilon_d>0$ depending only on the dimension $d$ so that
if $E$ has Hausdorff dimension at least ${d \over 2}$ then $D$ has Hausdorff dimension
at least ${1 \over 2} + \epsilon_d$.

\bigskip

To be precise, we prove Theorem 1.1, a multilinear version of Theorem 0.1. However, the interested
reader will see that Theorem 1.1 implies Theorem 0.1 following the treatment in [KT] (section 9).
In the case $d=2$, this theorem is the combination of the results of [KT] and [B], so that technically the
present paper is merely a generalization of sections 5 and 6 of [KT] to higher dimensions. This involves a few technical issues, but essentially uses the same ideas present there.

Largely, the inspiration for this paper is a preprint of Iosevich and Rudnev [IR], in which the authors
assert that a set having Hausdorff dimension ${d \over 2}$ but distance set of dimension ${1 \over 2}$
must have a lot of additive structure. This is obtained from a Fourier analytic estimate. We wrote this
paper because we felt that the authors understated the case. Not only does the set have additive
properties, but it also has multiplicative properties. Thus, {\it a fortiori}, apply Bourgain's sum-product
estimate, the set does not exist. It is the presence of this multiplicative structure in the higher dimensional case which the present work establishes

\head \S 1  Establishing linear structure \endhead

In this section, we will begin the proof of the ``multilinear
multidimensional Falconer conjecture."

Following [KT], we define the notion of a $(\delta, \rho)_d$ set.
We say that a set $A \subset {\Bbb R}^d $ is a $(\delta, \rho)_d$ set
provided that if $B$ is a ball of radius $r$, we have that
$$|B \cap A| \leq {\delta \over r}^{\rho}.$$ 

\bigskip

{\bf Theorem 1.1 (multilinear,multidimensional Falconer conjecture)} There is a constant $\epsilon_{d}>0$ so that
the following holds.  Let $B_1, \dots, B_{d+1}$ be unit balls 
in ${\Bbb R}^d$.
Let $E_1, \dots, E_{d+1}$ be $(\delta, {d \over 2})_d$ sets with
$E_j \subset B_j$. Suppose that for every $x_1,\dots,x_{d+1}$ with
$x_j \in B_j$ we have $|(x_2-x_1) \wedge (x_3-x_1) \dots \wedge 
(x_{d+1}-x_1)| \sim 1$. Let $D$ be a $(\delta,{1 \over 2})_1$ set.
Let
$$X=\{ (x_1,x_2, \dots , x_{d+1}) : x_j \in E_j, 1 \leq j \leq d+1,
|x_k-x_j| \in D, 1 \leq k < j \leq d+1 \}.$$
Then
$$|X| \leq  \delta^{ {d (d+1) \over 2} + \epsilon_{d} }.$$

\bigskip

We now say a few words about notation and the structure of the argument.
Our goal will be to show that a counterexample to the multilinear, multidimensional Falconer conjecture with a small value
of $\epsilon_d$ makes it possible to produce a counterexample to the
Bourgain sum-product theorem with $\epsilon=C_d \epsilon_d$ for some
dimensionally dependent constant $C_d$. We will make no effort to estimate this constant and will set up our notation accordingly. We shall denote by
$A \lessapprox B$, the expression $A \leq C \delta^{-\epsilon} B$, where
the constant and the $\epsilon$ vary from line to line. By subtle contrast,
we will use $A \lll B$ to mean that $A \leq C \delta^{\epsilon} B$ where the $C$ and $\epsilon$ vary from line but the $\epsilon$ is bigger than a large multiple of any $\epsilon$ that has appeared previously. A typical argument we will make will concern a set $S$ and a real valued function $q$ defined on $S$. If we know that 
$$|S| \gtrapprox 1,$$ 
and
$$|\{ x \in S: q(x) \lll 1\}| \lll 1, \tag 1.1$$
then
$$|\{x \in S: q(x) \gtrapprox 1\}| \gtrapprox 1. \tag 1.2$$
This is because the $\epsilon$ used in (1.1), while it must be large
compared to all previous $\epsilon$'s, need not be large compared
to the $\epsilon$ in (1.2). We use this idea to carry out various ``two ends" type arguments which show that of a certain class of configurations, most are nondegenerate. These same kinds of arguments are carried out repeatedly in [KT] by numbering constants and $\epsilon$'s. We find
the present notation simpler for both the author and the reader.

We begin with a geometrical lemma which help us to guarantee non-degeneracy.

\bigskip

{\bf Lemma 1.2 (Nondegeneracy)} With the notation of Theorem 1.1, let
$F \subset E_1$ and let $x_2 \in E_2, \dots, x_{d+1} \in E_{d+1}$.
Suppose that for every $y \in F$, we have that for each $2 \leq j \leq d+1$, we have $|y-x_j| \in D$. Then

(i) If $H$ is a $d-k$ dimensional affine subspace of ${\Bbb R}^d$. Then
$$|\{y \in F: \exists z \in H,  |z-y| \lessapprox \delta\}| \lessapprox
\delta^{{k + d\over 2}},$$
and
$$|\{y \in F: \exists z \in H,  |z-y| \lll 1\}| \lll
\delta^{{d \over 2}}.$$

(ii) If $S$ is a sphere whose center is a point $x \in E_2$ and
$H$ is a $d-k+1$ affine subspace of ${\Bbb R}^d$  then
$$|\{y \in F: \exists z \in H \cap S,  |z-y| \lessapprox \delta\}| \lessapprox
\delta^{{k + d\over 2}},$$
and
$$|\{y \in F: \exists z \in H,  |z-y| \lll 1\}| \lll
\delta^{{d \over 2}}.$$

\bigskip

{\bf Proof} To prove (i), we simply use
$$c(z)=(|z-x_2|,|z-x_3|, \dots, |z-x_{d+1}|)$$
as a coordinate system on $Q_1$. When restricted to $F$ these coordinates
take values only in the set $D$. By the implicit function theorem, everywhere locally on $H$ we may find a subset of cardinality $d-k$
of the coordinates which serve as a coordinate system.

To prove (ii), we do the same thing with
$$c(z)=(|z-x|,|z-x_3|,\dots,|z-x_{d+1}|).$$
Here the first coordinate is restricted to a single value and the others
vary in $D$. \qed

\bigskip

{\bf Remark 1.3} The same proof shows that $H \cap F$ is a
$({d-k \over 2}, \delta)_{d-k}$ set when viewed as a subset of $H$.
The same is lemma is true for any permutation of $1,\dots, k+1$
by the symmetry of the hypotheses. A set $F$ satisfying the hypotheses
of Lemma 1.2 will be referred to as {\bf nondegenerate}.

\bigskip

The following pigeonholing lemma will be used so frequently that is worth stating explicitly here:

\bigskip

{\bf Lemma 1.4 (Pigeonholing)} Let $A \subset {\Bbb R}^{n_1}$ and
$B \subset {\Bbb R}^{n_2}$ and let $\sim$ be a relation between them.
let $C=\{(a,b): a \in A,b \in B, a \sim b \}$. Then

i) {\bf (Chebychev)} There exists $b \in B$ so that
$$|\{ a \in A: a \sim b\}| \geq {|C| \over |B|}.$$

ii) {\bf (Cauchy-Schwarz)} 
$$|\{ (a,b,b^{\prime}) : a \sim b, a \sim b^{\prime}\}|
\geq {|C|^2 \over |A|}.$$

iii){\bf (H\"older)} 
$$|\{ (a,b_1, \dots, b_d) : a \sim b_1,b_2, \dots b_d\}|
\geq {|C|^d \over |A|^{d-1}}.$$

\bigskip

Of course Lemma 1.4 is proved   by applying respectively the
Chebychev, Cauchy-Schwarz, and H\'older inequalities.

Now we begin the proof of Theorem 1.1. We proceed by contradiction
and assume we are given $E_1, \dots, E_{d+1}$ which satisfy the hypotheses
of Theorem 1.1 but with
$$X=\{ (x_1,x_2, \dots , x_{d+1}) : x_j \in E_j, 1 \leq j \leq d+1,
|x_k-x_j| \in D, 1 \leq k < j \leq d+1 \},$$
we have
$$|X| \gtrapprox  \delta^{ {d (d+1) \over 2}  }.$$

\bigskip

{\bf Lemma 1.5(Refinement to nondegenerate sets)} We may find
$E_1^{\prime} \subset E_1$  and $E_2^{\prime} \subset E_2$ nondegenerate so that
$$Y=\{ (x_1,x_2, \dots , x_{d+1}) :x_1 \in E_1^{\prime}, 
x_2 \in E_2^{prime},  x_j \in E_j, 3 \leq j \leq d+1,
|x_k-x_j| \in D, 1 \leq k < j \leq d+1 \}.$$
Then
$$|Y| \gtrapprox  \delta^{ {d (d+1) \over 2} }.$$

{\bf Proof} It suffices to find $E_1^{\prime}$ nondegenerate so
that with
$$Z=\{ (x_1,x_2, \dots , x_{d+1}) :x_1 \in E_1^{\prime}, 
x_2 \in E_2,  x_j \in E_j, 3 \leq j \leq d+1,
|x_k-x_j| \in D, 1 \leq k < j \leq d+1 \},$$
we have
$$|Z| \gtrapprox  \delta^{ {d (d+1) \over 2} }.$$
Once we have done this, by symmetry, we can switch the subscripts 1 and
2 and apply the argument again.

We let
$$W=\{ (x_1,x_2, \dots , x_{d+1},y_2, \dots, y_{d+1}) : x_1 \in E_1$$
$$ x_j,y_j \in E_j, 2 \leq j \leq d+1,
|x_1 - x_j|, |x_1 - y_j||x_k-x_j|,|y_k-y_j| \in D, 2 \leq k < j \leq d+1 \},$$
Applying Lemma 1.4 (ii) we see that
$$|W| \gtrapprox \delta^{{(2d+1) d \over 2}}.$$
Now we apply Lemma 1.4 (i) to fix $y_2, \dots y_{d+1}$. Letting
$E_1^{\prime}$ be the set of $x \in E_1$ so that $|x-y_j| \in D$
for all $2 \leq j \leq d+1$, we obtain the desired conclusion. \qed

\bigskip

Our next few lemmas deal with the fact that the presence of a small distance
set implies that one can find many lines highly incident to our sets. The main
geometric fact to keep in mind is that the set of points in ${\Bbb R}^d$ which
are equidistant to a fixed set of $d$ points in general position $y_1, \dots, y_d$ is
the line perpindicular to the hyperplane containing the $y$'s which intersects that
hyperplane at the center of mass.

In the next series of lemmas for $t_1,t_2$ real numbers, we write $t_1 \sim t_2$ if
$$|t_1-t_2| \lessapprox \delta.$$

\bigskip

{\bf Lemma 1.6(Existence of lines)} Given $y_1, \dots , y_d$ in ${\Bbb R}^d$, we will say
they are $\in EGP$ (essentially in general position) provided that the $d-1$ vector
$(y_2-y_1) \wedge (y_3-y_1) \dots \wedge (y_d-y_1)$ is large, specifically
$$|(y_2-y_1) \wedge (y_3-y_1) \dots \wedge (y_d-y_1)| \gtrapprox 1.$$
Let
$$W=\{ (y_1,y_2, \dots , y_d,x): y_1,y_2 , \dots, y_d \in E_1^{\prime}, EGP, x \in E_2^{\prime},
|y_1-x| \sim |y_2-x| \sim \dots \sim |y_d-x|  \}.$$
Then
$$|W| \gtrapprox \delta^{{d^2+d -1 \over 2}}.$$

\bigskip

{\bf Proof} In light of Lemma 1.4 (i), we can fix $x_3, \dots, x_{d+1}$ in the conclusion of
Lemma 1.5, to see that
$$|\{ (x_1,x_2) \in E_1^{\prime} \times E_2^{\prime}  : |x_1-x_2| \in D \}| \gtrapprox \delta^{2d \over 2}. \tag 1.3$$
We define
$$E_2^{\prime \prime} = \{x \in E_2 : |\{x_1 \in E_1^{\prime}: |x_1 - x| \in D \}| \gtrapprox \delta^{d \over 2} \}.$$
By (1.3), we have
$$|E_2^{\prime \prime}| \gtrapprox |E_2^{\prime}| \gtrapprox \delta^{d \over 2}.$$
For a given $y  \in E_2^{\prime \prime}$ and a fixed $r \in D$, we have that
$$|\{x \in E_1^{\prime}: |x-y| \sim r \}| \lessapprox \delta^{d+1 \over 2},$$
by the nondegeneracy of $E_1^{\prime}$ (Lemma 1.2 (ii)).
Now since $D$ is a $(\delta, {1 \over 2})_1$ set, we must be able for each
$y \in E_2^{\prime \prime}$ to be able to find a set $D_{1} \subset D$ which is
$\ggg \delta$-separated and which has cardinality $\gtrapprox \delta^{-{1 \over 2}}$
so that for each $d \in D_1$ we have
$$|\{x \in E_1^{\prime}: |x-y| \sim d \}| \lessapprox \delta^{d+1 \over 2}.\tag 1.4$$
(We made $D_1$ to be $\ggg \delta$ separated so that the sets in (1.4) may be taken disjoint. That
is the $\sim$ in (1.4) still uses the $\epsilon$'s from before the $\ggg$, although the
lower bound on the cardinality of $D_1$ does not.).
Now we conclude that with
$$W_1=\{ (y_1,y_2, \dots , y_d,x): y_1,y_2 , \dots, y_d \in E_1^{\prime}, x \in E_2^{\prime \prime},
|y_1-x| \sim |y_2-x| \sim \dots \sim |y_d-x|  \},$$
we have
$$|W_1| \gtrapprox \delta^{{d^2+2d-1 \over 2}}.$$
This is not quite the proof of the lemma because $(y_1,\dots, y_d)$ are not necessarilly in 
essentially general position.
However we observe that by the nondegeneracy (Lemma 1.2 (ii) again) of $E^{\prime}$, we have for
$x \in E_2^{\prime \prime}$ and $r \in D$ fixed that
$$|\{(y_1,\dots,y_d) \in (E_1^{\prime}): |x-y_j| \sim r, |(y_2-y_1) \wedge \dots (y_d-1)| \lll 1
 \}| \lessapprox \delta^{d(d+1) \over 2}.\tag 1.5$$
 Thus is we define
$$W_2=\{ (y_1,y_2, \dots , y_d,x): y_1,y_2 , \dots, y_d \in E_1^{\prime},not EGP, x \in E_2^{\prime \prime},
|y_1-x| \sim |y_2-x| \sim \dots \sim |y_d-x|  \},$$
then we have
$$|W_2| \lll \delta^{{d^2+2d-1 \over 2}}.$$
Then we set
$$W=W_1 \backslash W_2,$$
so that we have
$$|W| \gtrapprox \delta^{{d^2+2d-1 \over 2}}.$$
\qed

We define $U$ to be the set of $d$-tuples $(y_1, \dots y_d) \in (E_1^{\prime})^d$ in $EGP$. Then
$$|U| \gtrapprox \delta^{{d^2 \over 2}}.$$ We say that $(y_1,\dots, y_d) \longrightarrow x$
(read $(y_1,\dots,y_d)$ point at $x$) if
$$|y_1-x| \sim |y_2-x| \dots \sim |y_d-x| \in D.$$

\bigskip

{\bf Lemma 1.7 (nondegeneracy of lines)}

(i) Fix $(y_1, \dots, y_d) \in U$. Then
$$|\{x:  (y_1,\dots,y_d) \longrightarrow x\}| \lessapprox \delta^{{2d-1 \over 2}}.$$

(ii) Fix $x \in E_2^{\prime \prime}$. Then
$$|\{(y_1,\dots,y_d) : (y_1,\dots,y_d) \longrightarrow x\}| \lessapprox \delta^{{d^2+d-1 \over 2}}.$$

\bigskip

{\bf Proof} To prove (i), we see that $\{x:  (y_1,\dots,y_d) \longrightarrow x\}$ is contained
in a $(\lessapprox \delta)$-tube. Thus we obtain the inequality immediately from nondegeneracy of
$E_2^{\prime}$.

To prove (ii), we observe that there is a collection of $\delta^{-{1 \over 2}}$ spheres $S_j$
centered at $x$
so that if $(y_1, \dots, y_d) \longrightarrow x$, then $y_1,\dots, y_d$ are all at distance
$\lessapprox \delta$ from $S_j$. Applying nondegeneracy of $E_1^{\prime}$ and summing
over $j$, we obtain the desired result.

\qed

\bigskip

We are in now a position to prove the analogue of the statement from [KT] which said
that the discretized Furstenburg conjecture implied the bilinear Falconer conjecture.

If $x_1,x_2 \in E_2^{\prime \prime}$ we say that $x_1 \triangle x_2$ (read $x_1$ has more than a
passing acquaintance with $x_2$) if there is $(y_1, \dots, y_d) \in U$ so that
$$(y_1, \dots, y_d) \longrightarrow x_1,x_2,$$
and
$$|\{x \in E_2^{\prime \prime}: (y_1, \dots, y_d) \longrightarrow x\}| \gtrapprox \delta^{{2d-1 \over 2}},$$
and
$$|x_1-x_2| \gtrapprox 1.$$

\bigskip

{\bf Lemma 1.8 (Furstenburg implies Falconer)}
We have the estimate
$$|\{ (x_1,x_2) \in (E_2^{\prime \prime})^2:  x_1 \triangle x_2 \}| \gtrapprox \delta^d.$$

\bigskip

{\bf Proof} We let $U^{\prime} \subset U$ be the set of $u \in U$ so that
$$|\{x \in E_2^{\prime \prime}: u \longrightarrow x\}| \gtrapprox \delta^{{d^2+d-1 \over 2}}.$$
In light of Lemma 1.6 and Lemma 1.7 (i), we see that
$$|U^{\prime}| \gtrapprox |U|,$$
and that
$$|\{ (x,u): u \in U^{\prime}, u \longrightarrow x\}| \gtrapprox \delta^{{d^2+2d-1 \over 2}}. \tag 1.5$$
Now we let $E_2^{(3)}$ be the set of $x$ in $E_2^{\prime \prime}$ so that
$$|\{u \in U^{\prime}: u \longrightarrow x \}| \gtrapprox \delta^{{d^2+d-1 \over 2}}.$$
In light of (1.5) and Lemma 1.7 (ii), we have that
$$|E_2^{(3)}| \gtrapprox |E_2^{\prime \prime}|.$$
We now demonstrate that points in $E_2^{(3)}$ have more than a passing acquaintance with
many points of $E_2^{\prime \prime}$.

Let $x \in E_2^{(3)}$. First we observe that for any $u$ with $u \longrightarrow x$, we have that
$$|\{x^{\prime} \in E_2^{\prime \prime} : u \longrightarrow x^{\prime}, |x-x^{\prime}| \lll 1\}|
\lll  \delta^{{2d-1 \over 2}}. \tag 1.6$$
If (1.6) failed and $u=(y_1,\dots, y_d)$ then $\{ |y_1-x^{\prime}|\} \subset D$ 
would fail to be a $(\delta, {1 \over 2})_1$ set where 
$x^{\prime}$ runs over the same set as in (1.6). Therefore for each $x \in E_2^{(3)}$ and
each $u \in U^{\prime}$ with $u \longrightarrow x$, we have
$$|\{x^{\prime} \in E_2^{\prime \prime}: |x-x^{\prime}| \gtrapprox 1, u \longrightarrow x^{\prime} \}|
\gtrapprox  \delta^{{2d-1 \over 2}}.\tag 1.7$$

Furthermore, for each $x \in E_2^{(3)}$, by Lemma 1.4 (i) there is a distance $r$ so that
$$|\{ u=(y_1,\dots,y_d) \in U^{\prime}: |x-y_1| \sim \dots \sim |x-y_d| \sim r, u \longrightarrow x\} |
\gtrapprox \delta^{{d^2+ d \over 2}}.$$
Again applying Lemma 1.4(i), we can fix $(y_2, \dots, y_d)$ so that
$$|\{ y_1 : u=(y_1,\dots,y_d) \in U^{\prime} , u \longrightarrow x\} |
\gtrapprox \delta^{{d+1 \over 2}}.$$
In this way, we obtain $\gtrapprox \delta^{{1-d \over 2}}$, choices of $y_1$ which are $\ggg \delta$
separated (here for grammatical reasons we have switched the order of $\ggg$ and $\gtrapprox$)
so that $u=(y_1, \dots, y_d) \in U^{\prime}$ and $u \longrightarrow x$. In light of (1.7), we see that
$$|\{x^{\prime} \in E_2^{\prime \prime}: x \triangle x^{\prime} \}| \gtrapprox \delta^{{d \over 2}}.$$
Now integrating over $x \in E_2^{(3)}$, we obtain the desired result. \qed

(Note that this last argument is a fairly standard synthesis of ``two ends" type conditions with the Bourgain bush argument.

From this point on we will essentially forget that $E_2^{\prime \prime}$ came from the
multilinear distance problem. We will use only the nondegeneracy of $E_2^{\prime \prime}$
together with the conclusion of Lemma 1.8.)

\head \S 2 Establishing additive and multiplicative structure \endhead

In this section, we begin with a set $E_2^{\prime \prime}$ which is nondegenerate and 
so that 
$$|\{ (x_1,x_2) \in (E_2^{\prime \prime})^2:  x_1 \triangle x_2 \}| \gtrapprox \delta^d. \tag 2.1$$

We begin by discretizing the set of lines implicit in the inequality (2.1). In what follows, we define
a $\delta$-tube to be the set of points at distance $\lessapprox \delta$ from a line. (This is a little
thicker than what is traditionally meant by a $\delta$-tube.)

\bigskip

{\bf Lemma 2.1} There is a set  $\Cal T$ of $\delta$-separated $\delta$-tubes with the property that
$$\#(\Cal T) \gtrapprox \delta^{1-d},$$
while
$$\sum_{T \in \Cal T} |T \cap E^{\prime \prime}| \gtrapprox \delta^{{1 \over 2}}.$$

\bigskip

{\bf Proof} We begin by stating an immediate consequence of nondegeneracy (Lemma 1.2)
which will be used frequently in this section.

\bigskip

{\bf Lemma 2.2} Let $T$ be a $\delta$-tube, then 
$$|T \cap E^{\prime \prime}|   \lessapprox \delta^{d-{1 \over 2}}.  $$

\bigskip

Next we observe that for any $x_1,x_2 \in E^{\prime \prime}$ with $x_1 \triangle x_2$,
there is a $\delta$-tube $T_{x_1 x_2}$ containing $x_1$ and $x_2$ so that
$$|T_{x_1 x_2} \cap E^{\prime \prime}| \gtrapprox \delta^{d-{1 \over 2}}.$$
We pick a maximal $\delta$-separated set of pairs $(x_1,x_2)$ and define
$\Cal T_0$ to be the associated tubes. By (2.1), we have that
$$\#(\Cal T_0) \gtrapprox \delta^{-d}.$$
Next we pick a maximal $\delta$-separated set of tubes in $\Cal T_0$ and refer to it as
$\Cal(T)$. By Lemma 2.2, we have that
$$\#(\Cal T) \gtrapprox \delta^{1-d}.$$

\qed

\bigskip

We observe that for any $x \in E^{\prime \prime}$, we have that
$$\#( \{ T \in \Cal T:  x \in T\} ) \lessapprox \delta^{{d-1 \over 2}}, \tag 2.2$$
since otherwise, we would have, by Bourgain's bush argument, that
$|E^{\prime \prime}| \ggg \delta^{{d \over 2}},$ which is impossible by nondegeneracy.

We define $F$ to be the set of $x \in E^{\prime \prime}$ so that
$$\#( \{ T \in \Cal T:  x \in T\} ) \gtrapprox \delta^{{d-1 \over 2}}.$$
We record an immediate consequence of (2.2), Lemma 2.1 and Chebychev's inequality :

\bigskip

{\bf Lemma 2.3}
$$|F| \gtrapprox |E^{\prime \prime}| \gtrapprox \delta^{{d \over 2}}.$$

\bigskip

We next observe that any refinement of $F$ will still satisfy the conclusion of Lemma 2.1

\bigskip
{\bf  Lemma 2.4}
Let $G \subset F$ with $|G| \gtrapprox |F|$ then
$$sum_{T \in \Cal T} |T \cap G| \gtrapprox \delta^{{1 \over 2}}.$$

\bigskip

{\bf Proof} Simply integrate $\#( \{ T \in \Cal T:  x \in T\} ) $ over $x$ in $G$.
\qed

\bigskip

We now find a refinement of $F$ which after a suitable projective transformation is contained
in a product set.

For $x,y \in F$, we define $x \square y$ provided $|x-y| \gtrapprox 1$ and there is 
$T \in \Cal T$ with $x,y \in T$.

We recall that $y_1, \dots, y_d$ are in essentially general position ($EGP$) if
$$|(y_2-y_1) \wedge \dots \wedge (y_d-y_1)| \gtrapprox 1.$$

\bigskip

{\bf Lemma 2.5} For $y_1, \dots , y_d \in F$, let
$$F(y_1, \dots, y_d) = \{x \in F: x \square y_1,y_2, \dots y_d\}.$$
Then
$$|\{(x,y_1, \dots, y_d): x \in F(y_1, \dots, y_d), y_1, \dots, y_d \in EGP \}| \gtrapprox \delta^{{d(d+1) \over 2}}.$$

\bigskip

{\bf Proof} We let $\tilde F$ be the set of $x$ in $F$ with 
$|\{y \in F: x \square y\}| \gtrapprox \delta^{{d \over 2}}.$
In light of nondegeneracy and Lemma 2.4, we have that $|\tilde F| \gtrapprox |F|$.
Thus we conclude
$$|\{ x,y \in F : x \square y\}| \gtrapprox \delta^{d}.$$
Applying H\"older's inequality (Lemma 1.4 (iii)), we get that
$$|\{x,y_1,\dots, y_d \in F: x \square y_1, \dots, y_d\}| \gtrapprox \delta^{{d(d+1) \over 2}}. \tag 2.3$$
Now applying nondegeneracy we get
$$|\{x,y_1,\dots, y_d \in F: x \square y_1, \dots, y_d, |(y_2-y_1) \wedge \dots \wedge (y_d-y_1)| \lll 1
 \}| \lll \delta^{{d(d+1) \over 2}}. \tag 2.4$$
 Comparing (2.3) with (2.4) gives the desired result.
 \qed
 \bigskip
 
 Applying Chebychev's inequality, we may find $y_1, \dots, y_d \in F$ in essentially
 general position so that $|F(y_1,\dots, y_d)| \gtrapprox |F|$. We let $F^{\prime}$ be the
 set of $x \in F(y_1,\dots, y_d)$ which are at distance $\gtrapprox 1$ from the hyperplane containing $y_1, \dots, y_d$. By nondegeneracy, we have that $|F^{\prime}| \gtrapprox |F|$.
 
 Next we find a suitable projective transformation. We imbed $\Bbb R^d$ in $\Bbb{RP}^d$
 by sending $(z_1, \dots, z_d)$ to $(z_1,\dots,z_d,1)$. We define $f_j$, the $j$th cardinal point
 at $\infty$ to be $(e_j,0)$ where $e_j$ is the $j$th standard basis vector in $\Bbb R^d$. We
 let $P$ be a (linear) projective transformation taking $y_j$ to $f_j$. Observe that $P$ is
 Lipschitz when restricted to $F^{\prime}$  with biLipschitz constant $\lessapprox 1$.
 
 We introduce some notation. Let 
 $$\pi_j: \Bbb R^d \longrightarrow \Bbb R,$$
 be the rank 1 projection which maps a point in $\Bbb R^d$ to its $j$th coordinate.
 Let
 $$\pi^j:\Bbb R^d \longrightarrow \Bbb R^{d-1},$$
 be the rank $d-1$ projection which maps a point in $\Bbb R^d$ to the point omits the $j$th
 coordinate.
 
 By construction, we know that for each $j$, we have
 $$|\pi^j ( P(F^{\prime}))| \lessapprox \delta^{{d-1 \over 2}},$$
 for each $j$. We would like to mimic the argument in [KT]. In order to do this we would need
 $$|\pi_j ( P(F^{\prime}))| \lessapprox \delta^{{1 \over 2}}.$$
 Fortunately, we may achieve this by picking a suitable refinement of $F^{\prime}$.
 (This is the only step which has no analog in the two dimensional case because $1+1=2$.
 
 \bigskip
 
 {\bf Lemma 2.6} Let $A \subset \Bbb R^d$ with 
 $$|A| \gtrapprox \delta^{{d \over 2}}.$$
 Suppose that for all $j$, we have that 
 $$|\pi^j (A)| \lessapprox \delta^{{d-1 \over 2}}.$$
 (Note that this immediately implies $|A| \lessapprox \delta^{{d \over 2}}$.)
 Then we may find $A^{\prime} \subset A$ with
 $$|A^{\prime}| \gtrapprox |A|,$$
 so that for each $j$, we have 
 $$|\pi_j(A)| \lessapprox \delta^{{1 \over 2}}.$$
 
 \bigskip
 
 {\bf Proof} We need only find a refinement $A^{\prime}$ with $|\pi_1 A| \lessapprox \delta^{{1 \over 2}}$
 because we can then permute the coordinates and repeat the procedure. Our main idea will be to
 observe that the map $\pi_1$ factors through the maps $\pi^2, \dots \pi^d$.
 We define $A_2 \subset A$  with $|A_2| \gtrapprox |A|$ and 
 so that for each $y \in \pi^2(A_2)$, we have $|\pi_2 ((\pi^2)^{-1}(y) \cap A       | \gtrapprox \delta^{{1 \over 2}}$.
 Then we define $A_3\subset A$  with $|A_3| \gtrapprox |A_2|$ and 
 so that for each $y \in \pi^3(A_2)$, we have $|\pi_3 (\pi^3)^{-1}(y) \cap A_2 | \gtrapprox \delta^{{1 \over 2}}$.
 We continue in the same fashion until we have defined $A_d=A^{\prime}$ so that
 for each $y \in \pi^d(A_d)$, we have $|\pi_d (\pi^d)^{-1}(y) \cap A_{d-1}| \gtrapprox \delta^{{1 \over 2}}.$
 Then since $\pi_1$ factors through each of $\pi^j$ for $j \geq 2$, by construction, we have
 for each $x \in A_d$ that 
 $$|\pi^1 (\pi_1^{-1} (\pi_1(x)) \cap A| \gtrapprox  \delta^{{d-1 \over 2}},$$
 so that by the upper bound on $|A|$,  we have
 $$|\pi_1(A_d)| \lessapprox \delta^{{1 \over 2}}.$$
 \qed
 
 \bigskip
 
 Applying Lemma 2.6 to $P(F^{\prime})$, mapping back under $P^{-1}$ and invoking the
 biLipschitz property of $P$, we see that we may choose $F^{\prime \prime} \subset
 F^{\prime}$ with 
 $$|F^{\prime \prime}| \gtrapprox |F^{\prime}|,$$
 so that
 $$P(F^{\prime \prime}) \subset A_1 \times A_2 \times \dots \times A_d,$$
 with each $A_j$ a $(\delta,{1 \over 2})_1$ set. (To obtain the $(\delta,{1 \over 2})_1$ property,
 observe that by the large amount of incidence with elements of $\Cal T$, there is
 a line whose intersection with $F^{\prime}$ is biLipschitz to refinements of each of $A_1,
 \dots, A_d$. Then we can invoke the nondegeneracy of $F^{\prime}$. We leave the
 details to the reader.)
 
 We are almost done. Our goal is now to show that $A_d$ has good additive and multiplicative
 properties, thereby contradicting Bourgain's sum product estimate. Our arguments now follow [KT]
 essentially verbatim.
 
 We record a version ({\it cf} [KT]) of the Gowers-Balog-Szemeredi theorem which we shall use.
 
 {\bf Lemma 2.7} Let $A,B \subset \Bbb R$ be bounded ($\lessapprox 1$) $(\delta,{1 \over 2})_1$
 sets with 
 $$|A|,|B| \gtrapprox \delta^{{1 \over 2}}.$$
 Suppose $G \subset A \times B$ with
 $$|G| \gtrapprox \delta,$$
 so that
 $$|\{ a+b: (a,b) \in G\}| \lessapprox \delta^{{1 \over 2}},$$
 then there exists $A^{\prime} \subset A$ with
 $$|A^{\prime}| \gtrapprox |A|,$$
 and
 $$|A^{\prime} - A^{\prime}| \lessapprox |A|.$$
 
 \bigskip
 
 We also record a version of Bourgain's theorem [B] which will provide us with our contradiction:
 
 {\bf Theorem 2.8(Bourgain)} Let $A \subset \Bbb R \backslash \{0\}$ with $A$ and ${1 \over A}$
 bounded and with $A$ a $(\delta,{1 \over 2})_1$ set satisfying
 $$|A| \lessapprox \delta^{{1 \over 2}}.$$
 Then either
 $$|A-A| \ggg \delta^{{1 \over 2}},$$
 or
 $$|AA| \ggg \delta^{{1 \over 2}}.$$
 
 \bigskip

We now recall where we are.

We have $(\delta,{1 \over 2})_1$ sets $A_1,\dots A_d$ so that
$P(F^{\prime \prime}) \cap A_1 \times \dots \times A_d$ is a refinement of $A_1 \times \dots \times A_d$.
We may pick a refinement $A^{\prime} \subset A_d$ with $|A^{\prime}| \gtrapprox |A_d|$
so that for any $A^{\prime \prime} \subset A^{\prime}$ with $|A^{\prime \prime}| \gtrapprox |A_d|$
we still have
$$|P(F^{\prime \prime}) \cap A_1 \times \dots \times A^{\prime \prime}| \gtrapprox
|A_1 \times \dots \times A_d|.$$
This will allow us to refine $A_d$ twice, once to ensure additve properties and once to ensure
multiplicative properties.

Most importantly, a large number of lines have large intersection with 
$A^{\prime} \times \dots \times A_d$. This is because of the projective and Lipschitz
properties of the map $P$. This if we apply $P$ to the part of $T$ at distance $\approx 1$
from the plane spanned by $y_1,\dots,y_d$, where $T \in \Cal T$, we get a subset of a
$\delta$-tube. This gives us a lot of arithmetic identities. 

To be precise, let $Z$ consist of the set of all pairs $(a_{11},\dots,a_{1d}),(a_{21},\dots,a_{2d})$
in $A_1 \times \dots \times A_d$ and all elements $a_{01} \in A_1$ with the
property that there exist $a_{02}, \dots, a_{0d}$ in $A_2,\dots,A_d$ respectively with
$$a_{1j}[{a_{01}-a_{11} \over a_{21} - a_{11}}] +a_{2j}[{a_{21}-a_{01} \over a_{21}-a_{11}}]  \sim
a_{0j},   \tag 2.5$$
for all $j$ from 2 to $d$. Then
$$|Z| \gtrapprox \delta^{{2d+1 \over 2}}.$$

We fix $j=d$. Then appropriately fixing all variables but $a_{1d}$ and $a_{2d}$
in (2.5) we see that
we can find 
$A^{\prime} \subset A_d$ so that $$|A_d -A_d| \lessapprox \delta^{{1 \over 2}}, \tag 2.6$$
using Lemma 2.7. Indeed we can do that same to $A_1$. We observe that we obtain many
equations of the form (2.5) running over only variables from the refined sets.

Then we fix all variables but $a_{1d}$ and $a_{01}$. We obtain many equations of the
form
$$C_1 a_{1d} a_{01} + C_2 a_{01} \sim a_{0d}.$$
Thus we may find refinements $A^{\prime \prime} \subset A^{\prime}$ so that
$$|(A^{\prime \prime } + {C_2 \over C_1})(A^{\prime \prime}+ {C_2 \over C_1})| \lessapprox \delta^{{1 \over 2}},$$
applying Lemma 2.7 multiplicatively (that is, to the logarithms)
However by (2.6), we see that 
$$|(A^{\prime \prime}  +{C_2 \over C_1}) + (A^{\prime \prime} + {C_2 \over C_1})|
\lessapprox \delta^{{1 \over 2}},$$
so that we have contradicted Theorem 2.8. Thus our assumption is false, and we have proven Theorem 1.1. \qed

\end